\documentclass[10pt]{amsart}
\usepackage{amsthm,amsmath,amssymb,amscd,graphics,enumerate, stmaryrd,xspace,verbatim, epic, eepic,url}
\usepackage[usenames,dvipsnames]{color}

\usepackage[colorlinks=true]{hyperref}

\usepackage[active]{srcltx}
\usepackage[all]{xypic}
\usepackage{mathrsfs}
\usepackage{tikz-cd}

\def\:{\colon}
\def\.{,\dots,}

\def\ZZ{\mathbf Z}
\def\RR{\mathbf R}
\def\NN{\mathbf N}
\def\AA{\mathbb A}

\def\ov#1{\overline{#1}}
\def\gp{\textrm{gp}}
\def\c#1{\mathcal{#1}}

\def\Supp#1{\textup{Supp}({#1})}
\def\Ker#1{\textup{Ker}(#1)}
\def\Span#1{\textup{Span}#1}

\newcommand{\Hom}{\operatorname{Hom}}
\newcommand{\Spec}{\operatorname{Spec}}
\newcommand{\invlim}{\displaystyle \lim_{ \longleftarrow } \,}
\newcommand{\dirlim}{\displaystyle \lim_{ \longrightarrow } \,}

{
   \newtheorem{theorem}[subsubsection]{Theorem}
      \newtheorem*{theorem*}{Theorem}

   \newtheorem{lemma}[subsubsection]{Lemma}

   \newtheorem{corollary}[subsubsection]{Corollary}
   
   \newtheorem*{conjecture*}{Conjecture}

}
{\theoremstyle{definition}
          \newtheorem*{exercise*}{Exercise}
   
   \newtheorem{example}[subsubsection]{Example}
   \newtheorem*{example*}{Example}
   
   \newtheorem{defn}[subsubsection]{Definition}
   \newtheorem*{definition*}{Definition}
   \newtheorem{rem}[subsubsection]{Remark}
   \newtheorem{remark}[subsubsection]{Remark}

}
\usepackage{hyperref}
\author{Sam Molcho}
\begin{document}
\title{Universal Stacky Semistable Reduction}
\maketitle

\bibliographystyle{amsalpha}

\begin{abstract}
Given a log smooth morphism $f: X \rightarrow S$ of toroidal embeddings, we perform a Raynaud-Gruson type operation on $f$ to make it flat and with reduced fibers. We do this by studying the geometry of the associated map of cone complexes $C(X) \rightarrow C(S)$. As a consequence, we show that the toroidal part of semistable reduction of Abramovich-Karu can be done in a canonical way.  
\end{abstract}

\begin{section}{Introduction}

The semistable reduction theorem of \cite{KKMS} is one of the foundations of the study of compactifications of moduli problems. Roughly, the main result of \cite{KKMS} is that given a flat family $X \rightarrow S = \Spec R$ over a discrete valuation ring, smooth over the generic point of $R$, there exists a finite base change $\Spec R' \rightarrow \Spec R$ and a modification $X'$ of the fiber product $X \times_{\Spec R} \Spec R'$ such that the central fiber of $X'$ is a divisor with normal crossings which is reduced. \\

Extensions of this result to the case where the base of $X \rightarrow S$ has higher dimension are explored in the work \cite{AK} of Abramovich and Karu. Over a higher dimensional base, semistable reduction in the strictest sense is not possible; nevertheless, the authors prove a version of the result, which they call weak semistable reduction. The strongest possible version of semistable reduction for a higher dimensional base was proven recently in \cite{ALT}. In all cases, the statement is proven in two steps. In the first, one reduces to the case where $X \rightarrow S$ is toroidal. This step requires resolution of singularities, and is proven in characteristic $0$. The second step is performing semistable reduction for the toroidal morphism; this is essentially a combinatorial problem, independent of characteristic. In \cite{AK}, this step has the following form: given a morphism $X \rightarrow S$ of projective toroidal embeddings, then one can find an alteration $S' \rightarrow S$ and a modification $X'$ of the normalization of the main component of $X \times_S S'$ which is weakly semistable, rather than semistable -- in other words, such that the family $X' \rightarrow S'$ is flat and has reduced fibers (and where $S'$ is non-singular). \\

The main result of this paper is that if we relax the hypotheses to allow families $X \rightarrow S$ of stacks rather than schemes, the toroidal part of weak semistable reduction can be done ``universally". Specifically, we show 
\begin{theorem}[Universal Weak Semistable Reduction]
\label{theorem:main}
Let $X \rightarrow S$ be any proper, surjective, log smooth morphism of toroidal embeddings. Then, there exists a commutative diagram 
\begin{align*}
\xymatrix{ \mathcal{X} \ar[r] \ar[d] & X \ar[d] \\ \mathcal{S} \ar[r] & S } 
\end{align*}
where $\mathcal{X} \rightarrow \mathcal{S}$ is a representable morphism of tame toroidal algebraic stacks, such that for any diagram
\begin{align*}
\xymatrix{ Y \ar[r] \ar[d] & X \ar[d] \\ T \ar[r] & S}
\end{align*}
where $T \rightarrow S$ a toroidal alteration; $Y$ is a modification of the normalization of the main component of the fiber product $X \times_S T$; and $Y \rightarrow T$ is weakly semistable, the morphism $Y \rightarrow T$ factors uniquely through $\mathcal{X} \rightarrow \mathcal{S}$. Furthermore, $\mathcal{X} \times_\mathcal{S} T \rightarrow T$ is weakly semistable.  
\end{theorem}  

The construction can be thought of as a logarithmic version of the Raynaud-Gruson flattening theorem. The stack structure ensures that the fibers of the morphism are reduced. However, it is not clear that the steps in our construction can be done separately -- first flattening the map without using stacks, then reducing the fibers. A similar idea of using stacks to perform semistable reduction universally in the case where the base is one dimensional was pursued by Olsson in \cite{Oss}, but our constructions are different.  

\end{section}

\begin{section}{The Toric Case}
We begin by studying the toric case first. We do this because the exposition is simpler in this case, yet all the essential ideas and proofs are already present. For simplicity, we work throughout over an algebraically closed field $k$. \\

Let $\mathcal{T}$ denote the category of toric varieties. We may identify a toric variety by a pair $(F,N)$ of a lattice $N$ and a fan $F$ in $N_\RR = N \otimes_\ZZ \RR$. We usually denote the toric variety associated to $(F,N)$ by $\AA(F,N)$, and refer to it as the geometric realization of $(F,N)$. We will blur the distinction between $(F,N)$ and $\AA(F,N)$ and refer to either as a toric variety depending on context. A morphism of toric varieties $(F,N) \rightarrow (G,Q)$ is a homomorphism of lattices $p:N \rightarrow Q$, such that $p_\RR:N_\RR \rightarrow Q_\RR$ takes each cone $\sigma \in F$ into a cone $\kappa \in G$. We recall the following notions, stated on the level of fans rather than on the level of varieties: 

\begin{defn}
The support $\Supp F$ of a fan $F$ is the set of vectors in $N_\RR$ that belong to some cone in $F$. 
\end{defn}   

\begin{defn}
A morphism of toric varieties $p:(F,N) \rightarrow (G,Q)$ is called proper if $p^{-1}(\Supp G) = \Supp F$. 
\end{defn}

\begin{defn}
A morphism of toric varieties $i:(F',N) \rightarrow (F,N)$ is called a modification or subdivision if $i:N \rightarrow N$ is the identity and $\Supp{F'} = \Supp F$.
\end{defn}

\begin{defn}
A morphism of toric varieties $j:(G',Q') \rightarrow (G,Q)$ is called an alteration if $j:Q' \rightarrow Q$ is a finite index injection and $\Supp {G'} = \Supp G$.
\end{defn}

\begin{rem}
Note that if $j:Q' \rightarrow Q$ is an injection with finite cokernel, then the condition $\Supp {G'} = \Supp G$ is equivalent to saying the morphism $(G',Q') \rightarrow (G,Q)$ is proper. Furthermore, for any homomorphism $j:Q' \rightarrow Q$ we get a 	``pull-back" fan $j_{\RR}^{-1}(G)$ which is isomorphic to $G$, since $j_\RR$ is an isomorphism. Thus, any toric alteration can be factored as a modification $(G',Q') \rightarrow (G,Q')$ composed with a finite index inclusion $(G,Q') \rightarrow (G,Q)$.  
\end{rem}

\begin{rem}
The pullback fan of remark 2.0.6  is a special case of the following more general construction: 
\begin{subsection}{Minimal Modification}
Let $p:N \rightarrow Q$ be a fixed homomorphism of lattices, and suppose $F,G$ are fans in the lattices $N,Q$ respectively. 

\begin{lemma}
\label{lem:minmod}
There exists a minimal modification $(F',N)$ of $(F,N)$ which maps to $(G,Q)$ inducing $p:N \rightarrow Q$.  
\end{lemma} 
\begin{proof}
The fan $F'$ is defined as $F':=\{p^{-1}(\kappa) \cap \sigma: \kappa \in G, \sigma \in F\}$. The proof that this is a fan and satisfies the universal property is straightforward and can be found in \cite{AM}. 
\end{proof}
\end{subsection}

\end{rem}

\begin{defn}
\label{def: weaklysemistable}
We call a morphism $p:(F,N) \rightarrow (G,Q)$ of toric varieties \emph{weakly semistable} if 
\begin{enumerate}
\item Every cone $\sigma$ in $F$ surjects onto a cone $\kappa \in G$ \\
\item Whenever we have $p(\sigma) = \kappa$, we have an equality of monoids $p(N \cap \sigma) = Q \cap \kappa$
\end{enumerate} 
\end{defn}

We remark here that our definition of weak semistability is slightly weaker than the one introduced in \cite{AK}: our base is not required to be smooth. Our definition is for practical purposes very close to the one of \cite{AK}. First, we have

\begin{lemma}
If $p:(F,N) \rightarrow (G,Q)$ is weakly semistable and $F,G$ are smooth, then $p$ is semistable. 
\end{lemma}

Next, it is proven in my appendix to \cite{Wmin} that the central fiber of such a morphism resembles the central fiber of a semistable morphism as in \cite{KKMS}, at least in codimension one. Furthermore, it is shown in \cite{AK} that a morphism that satisfies $(1)$ is equidimensional and a morphism that satisfies $(2)$ has reduced fibers. In \cite{AK}, equidimensionality is combined with smoothness of the base to deduce flatness of the morphism. Nevertheless, even over a singular base, we still have:

\begin{theorem}
\label{theorem:flat}
A morphism of toric varieties is weakly semistable if, and only if, it is flat and has reduced fibers.  
\end{theorem}

\begin{proof}
To begin, we will show that a weakly semistable morphism is saturated in the terminology of \cite{Ts}. From \cite[Theorem II.4.2]{Ts}, it will then follow that the map is flat with reduced fibers. The statement is local, so we may assume we are in the situation where a single cone $\sigma$ in $N$ maps into a single cone $\kappa$  in $Q$. By assumption, we have that faces of $\sigma$ map onto faces of $\kappa$, and whenever $\tau$ maps onto $\lambda$, we have $N \cap \tau$ mapping onto $\lambda \cap Q$. \\

Consider the dual monoids $Q_\kappa^{\vee}=\kappa^{\vee} \cap Q^{\vee}$ and $N_{\sigma}^{\vee} = \sigma^{\vee} \cap N^{\vee}$ in the dual lattices. Since $N \cap \sigma$ surjects onto $Q \cap \kappa$, the dual map $Q_\kappa^{\vee} \rightarrow N_{\sigma}^{\vee}$ is injective, and its cokernel has no torsion. To see flatness, we will verify that this dual map is an integral map of monoids in the sense of Kato. We use Kato's equational criterion for integrality \cite{K}. Suppose we are given  

\begin{align*}
p_1 + q_1 = p_2 + q_2
\end{align*} 

\noindent where $p_i \in N_\sigma^{\vee}$ and $q_i \in Q_\kappa^{\vee}$. We want to show that $p_1 = w + r_1$, $p_2 = w + r_2$, where $w \in N_\sigma^{\vee}$, $r_i \in Q_\kappa^{\vee}$, and $q_1+r_1 = q_2 + r_2$. Since the map $Q_\kappa^{\vee} \rightarrow N_{\sigma}^{\vee}$ is injective, we certainly have a (non-canonical) splitting of lattices $N_{\sigma}^{\vee,\gp} = Q_{\kappa}^{\vee,\gp}\oplus L$. So, we may identify any $p_1,p_2$ with $(w,r_1),(w,r_2)$ and we must have $q_1+r_1=q_2+r_2$. The point however is that this splitting may not respect the monoids, i.e $w$ may not be positive on $\sigma \cap N$. To fix this, we will carefully choose a particular splitting. Pick the face $\tau$ of $\sigma$ which maps isomorphically onto $\kappa$ and on which $p_1$ is minimal. To see that this is possible, let $v_1,\cdots,v_m$ be the extremal rays of $\kappa$, and let $u_k$ denote lifts of the rays $v_i$ in $\sigma$. Among the $u_k$, choose $u_1,u_2, \cdots u_m$ such that $u_i \mapsto v_i$ and such that $p_1(u_i)$ is minimal along all possible lifts of $v_i$ to an extremal ray of $\sigma$. The face $\tau$ of $\sigma$ generated by the $u_i$ is the desired face. By assumption, we have $\tau \cap N = Q \cap \kappa$. Using this splitting $N_\sigma^{\vee,\gp} = N_\tau^{\vee,\gp} \oplus L =  Q_\kappa^{\vee,\gp} \oplus L$, we see that we may write $p_1+q_1=p_2+q_2$ in the form $(w,r_1+q_1)=(w,r_2+q_2)$. We may identify every $\kappa$ with $\tau$, and thus the projection $\sigma \rightarrow \kappa$ gives us a map $p:\sigma \rightarrow \tau$. Every element $x$ of $\sigma$ can be written uniquely as $x=p(x)+v,v\in \Ker p$. Note that by construction, $w(x)=p_1(v)$, and $r_1(x)=p_1(p(x))$. To check that $w$ is non-negative on $\sigma$, it suffices to check it is non-negative on its extremal rays. For such a ray $x$ in $\sigma \cap \Ker p$, the result is clear since then $w(x) =p_1(x) \ge 0 $ by assumption. For an extremal ray not in $\Ker p$, we write $x = p(x)+v$, where $p(x)$ is an extremal ray on $\tau$. We then have that $p_1(p(x))\le p_1(x)$ by choice of $\tau$; hence $w(x) = p_1(v) = p_1(x-p(x)) \ge 0$, which completes the proof of integrality. The fact that the morphism of monoids is saturated then follows by \cite[Remark I.4.4]{Ts}.  

For the converse, suppose that the morphism $p: (F,N) \rightarrow (G,Q)$ is flat has reduced fibers. From \cite[Lemma 4.1]{AK}, it follows that all cones of $F$ map onto cones of $G$. Suppose then that the cone $\sigma$ in $N$ maps onto the cone $\kappa$ in $Q$, but that the map of lattices $N_\sigma:= N \cap \Span \sigma \rightarrow Q_\kappa = Q \cap \Span \kappa$ is not surjective. For a face $\tau < \sigma$ which maps isomorphically onto $\kappa$ (such a face must exist for dimension reasons), the map $N_\tau \rightarrow Q_\kappa$ is also not surjective, as the image of $N_\tau$ is contained in the image of $N_\sigma$. As the question is local, we may thus assume without loss of generality that $F$ consists of a single cone $\sigma$ which maps isomorphically to $\kappa$, and that $N_\sigma \rightarrow Q_\kappa$ therefore has finite cokernel. Choose a splitting $N = N_\sigma \oplus L$, $Q = Q_\kappa \oplus K$, and consider the diagram 
\begin{align*}
\xymatrix{\AA(\sigma, N) \cong \AA(\sigma,N_\sigma) \times \mathbb{G}_m^l \ar[r]^p \ar[d] & \AA(\kappa,Q) \cong \AA(\kappa,Q_\kappa) \times \mathbb{G}_m^l \ar[d] \\ \AA(\sigma, N_\sigma) \ar[r] & \AA(\kappa,Q_\kappa)}
\end{align*}
As $p$ and the two projections are flat with reduced fibers, the map $\AA(\sigma,N_\sigma) \rightarrow \AA(\kappa,Q_\kappa)$ must be flat with reduced fibers as well. So, as $\sigma \cong \kappa$ and $(N_\sigma)_\RR$ is isomorphic to $(Q_\kappa)_\RR$ we may further reduce to the case where $(F,N) \rightarrow (G,Q)$ is of the form $(\sigma, N) \rightarrow (\sigma, Q)$, with $\sigma$ full dimensional and $N \rightarrow Q$ a finite index inclusion, with non-trivial cokernel. But such a map is never reduced at the torus fixed point, and we obtain a contradiction. Therefore, $N_\sigma$ must surject onto $Q_\tau$, and $p$ is weakly semistable. 

\end{proof}

\begin{remark}
A consequence of \ref{theorem:flat} is that our definition of weak semistability, phrased in terms of the morphism of fans $(F,N) \rightarrow (G,Q)$, is equivalent to the definition of a saturated morphism of \cite{Ts}, a condition a priori phrased in terms of the log structures of $\AA(F,N),\AA(G,Q)$\footnote{This kind of translation between notions on a log scheme, which have more apparent functoriality properties, to notions on the fan/cone complex of the log scheme, which, at least to the author, provide better geometric intuition, is one of the themes of this paper}.  The key point of the proof is that weak semistability implies that the map of log structures is integral, which implies flatness, and then use condition $2$ in \ref{def: weaklysemistable} to show that the fibers are reduced. Note however that we are using both conditions of \ref{def: weaklysemistable} to prove integrality, not just condition $1$.   
\end{remark}

The advantage of this definition of weak semistability is that it is stable under pullbacks -- this is lemma \ref{lem:fiberproduct}. As the property of being flat with reduced fibers is clearly stable under pullback, in order to understand the meaning of the lemma, we need to discuss fiber products of toric varieties.  

\begin{subsection}{Fiber Products}\label{subsection: fiberproducts}

The category of toric varieties posseses fiber products: 

\begin{defn}[Toric Fiber Products]
The toric fiber product of 
\begin{align*}
\xymatrix{ & (F,N) \ar[d] \\ (H,L) \ar[r] & (G,Q) }
\end{align*}
is the toric variety with fan $F \times_G H = \{ \sigma \times_{\kappa} \lambda: \sigma \in F, \lambda \in H, \kappa \in G, \sigma \rightarrow \kappa,\lambda \rightarrow \kappa \}$ in the lattice $N \times_Q L $. 
\end{defn}
It is straightforward to verify that this collection forms a fan and that it satisfies the universal property of the fiber product with respect to toric maps. However, the toric fiber product can be ill-behaved, as it does not in general agree with the fiber product of the associated toric varieties in the category of schemes: 
\begin{align*}
\AA(F \times_G H) \ne \AA(F) \times_{\AA(G)} \AA(H)
\end{align*}

When considering the toric varieties as fine saturated log schemes, we can take yet another fiber product, the fiber product in the category of fine saturated log schemes, which we will denote by $(\AA(F) \times_{\AA(G)} \AA(H))_{\textup{tor}}$. This fiber product is closely related to the toric fiber product, though it is not exactly the same. We will explain the connection between $\AA(F \times_G H), \AA(F) \times_{\AA(G)} \AA(H)$ and $(\AA(F) \times_{\AA(G)} \AA(H))_\textup{tor}$ momentarily, after some preparation.  

\begin{example}
Consider the diagram 
\begin{align*}
\xymatrix{ & (\RR_+^2, \ZZ^2) \ar[d] \\ (\RR_+^2, \ZZ^2) \ar[r] & (\RR_+^2, \ZZ^2) }
\end{align*}
\noindent where the morphisms are $(a,b) \mapsto (a,a+b)$ and $(c,d) \mapsto (c+d,d)$ respectively (these are the two charts of the blowup of $\AA^2$ at the origin).  We have $\RR_+^2 \times_{\RR_+^2} \RR_+^2 = \{(a,b,c,d); a+b=d, a=c+d\}  \cong \RR_+ \subset \RR^2$. On the other hand, the fiber product of these two morphisms in the category of schemes is the variety $\{(x,y,z,w): xy = z, x = zw\}$ which is reducible.  
\end{example}

\begin{example}
What fails in the previous example is that the morphisms considered are not flat. However, flatness does not suffice to ensure toric fiber products agree with schematic fiber products. For instance, given the diagram 
\begin{align*}
\xymatrix{ & (\RR_+,\ZZ) \ar[d] \\ (\RR_+,\ZZ) \ar[r] & (\RR_+,\ZZ)  }
\end{align*}
where the two morphisms are $a \mapsto 2a$, $b \mapsto 3b$ respectively, the toric fiber product is $\{(a,b):2a=3b\} \cong \RR_+(3,2) \subset (\ZZ(3,2))_\RR$ whose geometric realization is $\AA^1$, whereas the schematic fiber product is $\{(x,y):x^2=y^3\}$. 
\end{example}

\begin{rem}[Colimits of Lattices]
Given a diagram of lattices, we may take the limit or colimit of the diagram in the category of abelian groups. The limit of such a diagram is always a lattice, hence coincides with the limit in the category of lattices as well. In general, colimits of lattices are not lattices. However, given a finitely generated abelian group $L$, we can form the associated lattice $L/\textup{Torsion}$. The functor $L \mapsto L/\textup{Torsion}$ is a left adjoint, and thus, the colimit of the $L_i$ in the category of lattices coincides with $\dirlim L_i/\textup{Torsion}$, where the colimit is understood in the category of abelian groups. 
\end{rem}

\begin{rem}[Double Dual of a Lattice]
Note that for any finitely generated abelian group $L$, we have a natural isomorphism $L/\textup{Torsion} \cong (L^{\vee})^{\vee}$, where $L^{\vee} = \Hom(L,\ZZ)$ as usual. Since we have 
\begin{align*}
(\dirlim L_i)^{\vee}:=\Hom(\dirlim L_i,\ZZ) = \invlim L_i^{\vee} 
\end{align*}
by the defining property of a colimit, it follows that 
\begin{align*}
(\invlim L_i)^{\vee} = (\dirlim L_i^{\vee})^{\vee})^{\vee} \cong \dirlim L_i^{\vee}/\textup{Torsion}
\end{align*}
\noindent In particular, we have that for lattices $Q,N,L$ 
\begin{align*}
(N \times_Q L)^{\vee} = N^{\vee} \oplus_{Q^{\vee}} L^{\vee}/\textup{Torsion}
\end{align*}
In the situations we are interested, the map $N \rightarrow Q$ will arise from a log smooth morphism. Therefore, the map $Q^{\vee} \rightarrow N^{\vee}$ will be injective, and in characteristic $p$ its cokernel will not have $p$-torsion. We claim that in this situation, the pushout has no $p$-torsion either. This is essentially because log smoothness is stable under pullbacks, but one can see it directly by applying the functor $\Hom(\ZZ/p\ZZ,?)$ to the exact sequence
\begin{align*}
\begin{xymatrix}
{0 \ar[r] & L^{\vee} \ar[r] & N^{\vee} \oplus_{Q^{\vee}} L^{\vee} \ar[r] & N^{\vee}/Q^{\vee} \ar[r] & 0} 
\end{xymatrix}
\end{align*} \\
\end{rem}
In general, for cones $\sigma \in F, \kappa \in G, \lambda \in H$, we get maps $\sigma^{\vee} \cap N^{\vee} = N_{\sigma}^{\vee} \rightarrow N^{\vee}$, $Q_{\kappa}^{\vee} \rightarrow Q^{\vee}, L_{\lambda}^{\vee} \rightarrow L^{\vee}$ and so a map $N_{\sigma}^{\vee} \oplus_{Q_{\kappa}^{\vee}} L_{\lambda}^{\vee} \rightarrow N^{\vee} \oplus_{Q^{\vee}}L^{\vee} \rightarrow (N \times_Q L)^{\vee} = N^{\vee} \oplus_{Q^{\vee}}L^{\vee}/\textup{Torsion}$. It is clear that vectors in the image of $N_{\sigma}^{\vee} \oplus_{Q_{\kappa}^{\vee}} L_{\lambda}^{\vee}$ are non-negative on $\sigma \times_{\kappa} \lambda \subset N \times_Q L$, and so we get a map $N_{\sigma}^{\vee} \oplus_{Q_{\kappa}^{\vee}} L_{\lambda}^{\vee} \rightarrow (\sigma \times_{\kappa} \lambda)^{\vee}$. On the one hand, 
\begin{align*}
k[(\sigma \times_{\kappa} \lambda)^{\vee} \cap (N \times_Q L)^{\vee}]
\end{align*}
are the affine charts for the fiber product of $(F,N),(G,Q),(H,L)$ in the category of toric varieties. On the other hand, since the functor \textbf{Mon} $\rightarrow$ \textbf{$k$-alg}, $M \rightarrow k[M]$ is left adjoint to the inclusion of $k$-algebras into monoids, the functor preserves colimits, so 
\begin{align*}
k[N_{\sigma}^{\vee} \oplus_{Q_{\kappa}^{\vee}} L_{\lambda}^{\vee}] \cong k[N_{\sigma}^{\vee}] \otimes_{k[Q_{\kappa}^{\vee}]} k[L_{\lambda}^{\vee}]
\end{align*}
\noindent which are the affine charts of the fiber product in the category of schemes. Thus, we see that the toric fiber product and the usual schematic fiber product coincide if and only if 
\begin{align*}
N_{\sigma}^{\vee} \oplus_{Q_{\kappa}^{\vee}} L_{\lambda}^{\vee} \rightarrow (\sigma \times_{\kappa} \lambda)^{\vee} \cap (N \times_Q L)^{\vee}
\end{align*} 
is an isomorphism for all cones $\{\sigma \times_{\kappa} \lambda\}$ in $F \times_G H$. The monoid $(\sigma \times_{\kappa} \lambda)^{\vee} \cap (N \times_Q L)^{\vee}$ can be obtained from $N^{\vee} \oplus_{Q^{\vee}}L^{\vee}$ as a composition of three functors: It is the composition of the integralization functor $P \mapsto P^{\textup{int}}$ which replaces $P$ with its image in $P^\gp$, followed by the saturation functor $P \mapsto P^{\textup{sat}} = \{x \in P^\gp: nx \in P^{\textup{int}}\}$, followed by the functor $P \mapsto P/\textup{Torsion}$. As $P \rightarrow P^{\textup{int}}$ is surjective, $\Spec k[P^{\textup{int}}] \rightarrow \Spec {k[P]}$ is a closed immersion, and as $\Spec k[P^\gp]$ is open in $\Spec k[P^{\textup{int}}]$, the geometric realization $\Spec k[P^{\textup{int}}]$ of $P^{\textup{int}}$ is the closure of $\Spec k[P^\gp]$ in $\Spec k[P]$. Similarily, given an integral monoid $P \subset P^\gp$, the algebra $k[P^{\textup{sat}}] \subset k[P^\gp]$ is the integral closure of $k[P]$, and hence $\Spec k[P^{\textup{sat}}] $ is the normalization of $\Spec k[P]$. Performing these two functors produces the fiber product in the category of f.s. log schemes. Therefore, the fiber product $(\AA(F) \times_{\AA(G)} \AA(H))_\textup{tor}$ is the normalization of the main component of $\AA(F) \times_{\AA(G)} \AA(H)$, that is, of the component containing the fiber product of the tori. It is also almost the toric fiber product $\AA(F \times_G H)$, but it is not quite as $\Spec k[P^\gp]$ is not a torus if there is torsion. However, in the situation we are interested in, where the torsion is prime to the characteristic of $k$, the process of killing torsion is very mild, and exhibits $\Spec k[P]$ as a disjoint union of schemes isomorphic to $\Spec k[P/\textup{Torsion}]$. Therefore, the toric fiber product $\AA(F \times_G H)$ is obtained from the log fiber product $(\AA(F) \times_{\AA(G)} \AA(H))_{\textup{tor}}$ by keeping only the component that contains the identity of the fiber product of tori. \\

With these observations at hand, we find:  

\begin{lemma}
\label{lem:fiberproduct}
Suppose 
\begin{align*}
\xymatrix {(F,N) \ar[d]_p \\ (G,Q)}
\end{align*}
is weakly semistable. Then, for a map $(G',Q') \rightarrow (G,Q)$ of toric varieties, the geometric realization of the diagram 
\begin{align*}
\xymatrix{ F_\tau = F \times_G G' \ar[r] \ar[d]_{p_G'} & F \ar[d]^p \\ G' \ar[r] &  G}
\end{align*}
is cartesian in the category of schemes, and $p_{\tau}$ is also weakly semistable. 
\end{lemma}

\begin{proof}
The property of being flat with reduced fibers is stable under pullbacks. Thus, the interesting point is showing that the diagram is cartesian in the category of schemes. Both properties are local on $F$, so we may replace $F$ by a single cone $\sigma$, $G$ by a single cone $\kappa$, and $G'$ by a single cone $\lambda$.
Since $p$ is flat with reduced fibers, theorem II.4.2 in \cite{Ts} implies that the pushout 
\begin{align*}
N_{\sigma}^{\vee} \oplus_{Q_{\kappa}^{\vee}} L_{\lambda}^{\vee}
\end{align*}
is saturated in its associated group. This associated group is actually a lattice, since by choosing a face of $\tau$ of $\sigma$ which maps isomorphically to $\kappa$, and with $N \cap \tau \cong Q \cap \kappa$, we obtain a splitting $(N_{\sigma}^{\vee})^\gp \cong (Q_{\kappa}^{\vee})^\gp \oplus N'$ -- c.f  the proof of \ref{theorem:flat}. Thus, 
\begin{align*}
N_{\sigma}^{\vee} \oplus_{Q_{\kappa}^{\vee}} L_{\lambda}^{\vee}
\end{align*}
is identified with the intersection of a cone $C$ in $(N \times_Q L)_{\RR}^{\vee}$ with $(N \times_Q L)^{\vee}$. The same is true with $(\sigma \cap \lambda) \cap (N \times_Q L)$. However, the dual of each of these cones is isomorphic to $\sigma \times_{\kappa} \lambda$, in the first case by the defining property of the colimit, and in the second by the relation $(C^{\vee})^{\vee}=C$ for the double dual of a cone in a fixed lattice. Applying the dual again, we see that $N_{\sigma}^{\vee} \oplus_{Q_{\kappa}^{\vee}} L_{\lambda}^{\vee}$ must be isomorphic to $(\sigma \cap \lambda) \cap (N \times_Q L)$, and the result follows by the discussion preceeding the lemma. 
\end{proof}

\end{subsection}
\begin{subsection}{Toric Stacks}
In what follows, we will need the notion of a \emph{toric stack}. For us, a toric stack will be always given by the data of a ``KM" fan, i.e a triple $(F,N,\{N_{\sigma}\}_{\sigma \in F})$, where $(F,N)$ is the usual data of a toric variety, and $N_{\sigma}$ is a collection of sublattices of $N$, one for each $\sigma \in F$, with the properties 
\begin{itemize}
\item $N_{\sigma} \subset N \cap \Span(\sigma)$ is a finite index inclusion. \\
\item $N_{\sigma} \cap \Span \tau = N_{\tau}$ for a face $\tau$ of $\sigma$. 
\end{itemize}
A morphism of KM fans $(F,N,\{N_{\sigma}\}) \rightarrow (G,Q,\{Q_{\kappa}\})$ is a morphism $(F,N) \rightarrow (G,Q)$ such that whenever $\sigma \mapsto \kappa$, $N_{\sigma} \rightarrow N \rightarrow Q$ factors through $Q_{\kappa}$.  \\
 
The data of a KM fan $(F,N,N_{\sigma})$ has a geometric realization into a stack $\AA(F,N,\{N_{\sigma}\})$. The stack is constructed as follows: the finite index sublattice $N_\sigma \rightarrow N \cap \Span \sigma$ lifts to a finite index sublattice $L \rightarrow N$ with $L \cap \Span{\sigma} = N_\sigma$, which induces a finite map $T(L) \rightarrow T(N)$ of tori, with kernel a finite group $K_\sigma$. Then $\AA(F,N,\{N_{\sigma}\})$ is the colimit of the local pieces $\AA(\sigma,N,N_\sigma):=[\AA(\sigma,L)/K_\sigma]$\footnote{this quotient is independent of the choice of the lift $L$}, with a map $\AA(\tau,N,N_\tau) \rightarrow \AA(\sigma, N, N_\sigma)$ for every face $\tau <\sigma$. The coarse moduli space of the stack is the toric variety $(F,N)$. The stack is a tame Artin stack, and is Deligne-Mumford whenever the index of $N_\sigma$ in $N$ is prime to the characteristic for all cones $\sigma$ -- so in particular, this is always the case in characteristic $0$. The data of a morphism of KM fans has a geometric realization into a morphism of the geometric realizations of the fans, and the associated morphism between coarse spaces is simply the geometric realization of $(F,N) \rightarrow (G,Q)$. The notion of a KM fan was first introduced in \cite{Ttor}; their main properties and the connection of the geometric realization with intrinsically defined toric stacks are developed in \cite{GMtor}. In the language of log geometry, the data of a stacky fan determines a Kummer extension of the log structure of $\AA(F,N)$, and $\AA(F,N,\{N_\sigma\})$ is exactly the stack one obtains from the construction of \cite{BV} -- we will use this observation later to define toroidal stacks, but we will phrase everything using the fan formalism as it is better suited for doing geometry on the combinatorial level.  \\

Any toric variety can be regarded as a toric stack, by taking $N_{\sigma} = N \cap \Span \sigma$ for each cone $\sigma$ -- note that there is no additional information in the $N_{\sigma}$ in this case. Under this identification, the category of toric varieties becomes a full subcategory of the category of toric stacks. We will use this identification in what follows and keep denoting a toric variety $(F,N)$ by $(F,N)$ even when the context makes it clear that it is considered as a toric stack. The proof of the following lemma is simple, and can be found in \cite{GMtor}:  
\begin{lemma}
\label{lem:representable}
Let $p:(F,N,\{N_{\sigma}\}) \rightarrow (G,Q,\{Q_{\kappa}\})$ be a morphism of KM fans. If $p^{-1}(Q_{\kappa}) = N_{\sigma}$ whenever $\sigma \mapsto \kappa$, the geometric realization of $p$ is representable.  
\end{lemma}
\end{subsection}

\begin{subsection}{The Main Construction}
We now fix the morphism $p:(F,N) \rightarrow (G,Q)$ which is \emph{surjective and proper}. 
\begin{defn}
\label{defn:main}
Let $\mathcal{C}$ be the category whose objects are diagrams
\begin{align*}
 \xymatrix{(\Phi,N') \ar[r]^j \ar[d]_\pi & (F,N) \ar[d]^p \\ (\Gamma,Q') \ar[r]_i & (G,Q)}
\end{align*}
such that 
\begin{itemize}
\item The map $i$ is an alteration. \\
\item $N'$ is the fiber product $N \times_Q Q'$. \\ 
\item $\Phi$ is a modification of $j^{-1}(F).$\\
\item $\pi$ is weakly semistable. \\
\end{itemize}
A morphism in $\mathcal{C}$ is a commutative diagram 
\begin{align*}
\xymatrix{(\Phi'',N'') \ar[r] \ar[d] & (\Phi',N') \ar[d] \\ (\Gamma'',Q'') \ar[r] & (\Gamma',Q')}
\end{align*}
\noindent which commutes with the morphisms to $p:(F,N) \rightarrow (G,Q)$.
\end{defn} 

\begin{theorem} 
\label{theorem:maintoric}
The category $\mathcal{C}$ has a terminal object which is a toric stack. In other words, there is a diagram 
\begin{align*}
\xymatrix{(F',N,\{N_\sigma\}) \ar[r] \ar[d] & (F,N) \ar[d] \\ (G',Q,\{Q_\kappa\}) \ar[r] & (G,Q)}
\end{align*}
such that every diagram
\begin{align*}
\xymatrix{(\Phi,N') \ar[r]^j \ar[d]_\pi & (F,N) \ar[d]^p \\ (\Gamma,Q') \ar[r]_i & (G,Q)}
\end{align*}
factors uniquely as 
\begin{align*}
\xymatrix{(\Phi,N') \ar[r] \ar[d]_\pi & (F',N,\{N_\sigma\}) \ar[r] \ar[d] & (F,N) \ar[d]^p \\ (\Gamma,Q') \ar[r]  & (G',Q,\{Q_\kappa\})  \ar[r]& (G,Q)}
\end{align*}
\end{theorem}

\begin{proof}
We first construct $G'$. The idea of the construction comes from \cite{KSZ}. Let $p(F)$ denote the collection of images of cones of $F$. Note that though every cone $p(\sigma)$ is contained in a cone of $G$, thus is convex, $p(F)$ is in general not a fan, as cones may not intersect along faces. We define $G'$ as the subdivision of $G$ determined by the cones in $p(F)$. Explicitly, this means the following: For every vector $w$ in $G$, we look at the collection 
\begin{align*}
N_0(w) = \{\sigma \in F: p(v)=w \textup{ for some } v \textup{ in the interior } of \sigma\}
\end{align*}
The cones $\kappa$ of $G'$ are precisely the cones such that for any two $w,w'$ in the interior of $\kappa$, we have $N_0(w)=N_0(w')$. Next, for a cone $\kappa \in G'$, we take 
\begin{align*}
Q_\kappa = \cap_{\sigma \in N_0(\kappa)} p(N \cap \sigma)
\end{align*}
In the interior of $\kappa$, this has the following description: $w \in Q_\kappa$ if and only if there exist $v_i \in \sigma_i$ with $p(v_i)=w$ for every cone $\sigma_i \in N_0(\kappa)$.  This completes the construction of the base $(G',Q,Q_\kappa)$. \\

At this point we need to verify that this construction actually yields a fan. The difficult part is verifying that the cones are strictly convex. So fix a cone $\kappa \in G'$, and pick two interior vectors $w,w' \in \kappa$. We will show that the whole line segment connecting $w$ to $w'$ must also be in the interior of $\kappa$. Suppose there exists a $t \in (0,1)$ for which $N_0(tw+(1-t)w')$ is different from $N_0(w)=N_0(w')$. Take for simplicity the smallest such $t$ -- this makes sense since the condition $N_0(w)=N_0(u)$ is an open condition on $u$-- and denote the point $tw+(1-t)w'$ by $w''$ to ease the notation. Certainly, since every cone $\sigma_i$ in $N_0(w)$ is strictly convex, the line segment between two lifts of $w,w'$ in $\sigma_i$ is also in $\sigma_i$, so $N_0(w) \subset N_0(w'')$. So take a cone $\sigma \in N_0(w'')-N_0(w)$, and a lift $v''$ of $w''$ in $N_0(w'')$. We look at the fiber of the map of vector space $N_{\RR} \rightarrow Q_{\RR}$ over the interval $[w,w']$ in $Q_{\RR}$. Call this fiber $N_{[w,w']}$, and let   $F_{[w,w']}$ be the intersection of $N_{[w,w']}$ with the fan $F$. Then $F_{[w,w']}$ is a polyhedral decomposition of $N_{[w,w']}$. The cone $\sigma$ intersects $F_{[w,w']}$ into a one-dimensional union of cells, since the relative dimension of $\sigma$ under $F \rightarrow G$ is $0$ by assumption, and contains $v''$ as an extreme point. Similarily, cones in $N_0(w)$ correspond to edges in $F_{[w,w']}$. Since $(F,N) \rightarrow (G,Q)$ is surjective and proper, the support of $F_{[w,w']}$ is all of $N_{[w,w']}$. In particular, the star of $v''$ in $F_{[w,w']}$ must intersect $F_{[w,w'')}$ non-trivially; so, in particular, there is an edge in the star of $v''$ in $F_{[w,w']}$ which maps to a vector $sw+(1-s)w'$ with $s<t$. By assumption on $t$, we have $N_0(sw+(1-s)w')=N_0(w)$, so in fact the edge corresponds to a cone in $N_0(w)$ and thus contains a lift of $w$. Since $\sigma$ is by choice not in $N_0(w)$, $v''$ is an extreme point of the edge. But this is a contradiction, since the edge must extend to contain a lift of $w'$ as well, as we assumed $N_0(w)=N_0(w')$. Thus we must have $N_0(w)=N_0(sw+(1-s)w')$ for all $s \in [0,1]$, and convexity follows.  

To construct $(F',N,\{N_\sigma\})$, we simply take the minimal subdivision of $(F,N)$ that maps to $(G',Q,Q_\kappa)$, as in \ref{lem:minmod}.  This means that $F'$ is the fan $\{p^{-1}(\kappa) \cap \sigma: \kappa \in G', \sigma \in F\}$ and the sublattice corresponding to $\sigma' := p^{-1}(\kappa) \cap \sigma$ is $L_{\sigma'} := p^{-1}(Q_\kappa \cap \sigma')$. Furthermore, the morphism $(F',N,\{N_\sigma\})$ is weakly semistable: cones of $F'$ map onto cones of $G'$ by construction, and $N_{\sigma}$ maps onto $Q_{\kappa}$ whenever $\sigma \mapsto \kappa$ by construction again. \\

Suppose now we are given a diagram 

\begin{align*}
\xymatrix{(\Phi,N') \ar[r]^j \ar[d]_\pi & (F,N) \ar[d]^p \\ (\Gamma,Q') \ar[r]_i & (G,Q)}
\end{align*}

\noindent where $i$ is an alteration, $N'$ the fiber product $N \times_Q Q'$, and $\Phi$ a subdivision of $j^{-1}F$. Assume furthermore that $\pi$ is semistable. Let $w,w'$ be two lattice points in the interior of a cone $\gamma$ of the fan $\Gamma$. Suppose that $w$ maps into a cone $\kappa \in G'$; we show that $w'$ maps to the same cone as well. Consider lifts $v_1,\cdots,v_n$ of $v=i(w)$ to cones $\sigma_i \in N_0(w) \subset F$. Since $\Phi$ subdivides $j^{-1}F$, there are cones $g_1,\cdots,g_n$ in $\Phi$ such that $j(g_i) \subseteq \sigma_i$; so we may find lifts $w_1,\cdots,w_n$ of $w$ in $g_i$. But then each cone $g_i$ maps to $\gamma$ under the projection $\pi$, and hence maps onto $\gamma$ and $\gamma \cap Q' = \pi(g_i \cap N)$ from conditions $(1),(2)$ in the definition of semistability. Since $w'$ is in $\gamma \cap Q'$ as well, this means that there exists $w_1',\cdots,w_n' \in g_i \cap N$ that map to $w'$ as well -- and hence there are $v_1'=j(w_1'),\cdots,v_n'=j(w_n')$ in the cones $\sigma_i \in N_0(w)$ that map to $w'$ as well. It follows that $N_0(w) \subset N_0(w')$, thus, by symmetry, $N_0(w)=N_0(w')$; hence $w,w'$ belong to the same cone $\kappa$ of $G'$. Furthermore, they are in the image of the lattice $N_{\sigma_i}$ for each cone in $N_0(w)$, thus in fact in the monoid $Q_\kappa$. Thus $(\Gamma,Q')$ factors through $(G',Q,Q_\kappa)$. The fact that $(\Phi,N)$ must factor through $(F,N,N_\sigma)$ factors through the universal property defining $(F,N,N_\sigma)$ automatically.

\end{proof}

\begin{rem}
It is worth pointing out that this proof goes through without assuming that the map $i:(\Gamma,Q') \rightarrow (G,Q)$ is an alteration. All that is required is that the kernel of $N \rightarrow Q$ and that the kernel of $N' \rightarrow Q'$ coincide. 
\end{rem}

Using the notation of definition \ref{defn:main}, we have as a corollary: 
\begin{corollary}
The minimal modification $\Phi$ of $j^{-1}(F)$ such that
\begin{align*}
 \xymatrix{(\Phi,N') \ar[r]^j \ar[d]_\pi & (F,N) \ar[d]^p \\ (\Gamma,Q') \ar[r]_i & (G,Q)}
\end{align*}
\noindent commutes and $\pi$ is weakly semistable is given by the fiber product of 
\begin{align*}
\xymatrix{& (F,N,\{N_{\sigma}\}) \ar[d]^p \\ (\Gamma,Q') \ar[r]_i & (G',Q,\{Q_{\kappa}\})}
\end{align*}
\noindent Its geometric realization coincides with $\AA(\Gamma,Q') \times_{\AA(G',Q,\{Q_{\kappa}\})} \AA(F,N,\{N_{\sigma}\})$.
\end{corollary}
\end{subsection}

\begin{proof}
This follows immediately by combining \ref{theorem:maintoric}, \ref{lem:fiberproduct}, \ref{lem:representable}. 
\end{proof}

\end{section}

\begin{section}{Globalizing}
We are now ready to discuss the changes necessary to generalize the above construction to the toroidal case. Recall the relevant definitions from \cite{KKMS}. To any toroidal embedding $(X,U)$ there is associated a stratification, the strata being determined by the irreducible components of the divisor $X-U$. For each stratum $Y$, we set
\begin{align*}
M^Y = \textup{Divisors on } \textup{Star} Y \\
M^Y_+ = \textup{Effective Divisors on } \textup{Star} Y \\
N^Y = \Hom(M_Y,\ZZ)\\
\sigma^Y = \{v \in N^Y \otimes_\ZZ \RR: v \textup{ is non-negative on } M^Y_+\}
\end{align*}
The collection of the cones $\sigma^Y$ together with their integral structure $N^Y$ is a cone complex, which we denote by $C(X)$; the only contrast with the toric theory is that the cones do not all inhabit a single (canonical) lattice $N$. Subdivisions of this cone complex correspond to birational modifications of $(X,U)$ which are the identity on $U$. \\

At this point, it is useful to connect the theory of toroidal embeddings with the language of logarithmic geometry\footnote{We will not use logarithmic geometry in a deep way. However, to avoid the language completely would make the presentation of what follows unnecessarily cumbersome. The reader who wishes to circumvent the use of logarithmic geometry entirely should be able to prove all relevant statements by using analogous results in \cite{KKMS} without much difficulty.}. A toroidal embedding $(X,U)$ carries a canonical structure of a log scheme $(X,M_X)$, by setting 
\begin{align*}
M_X(V) = \{f \in \mathcal{O}_X(V): f \in \mathcal{O}^*_X(V \cap U)\}
\end{align*}
Conversely, by the chart criterion of log smoothness of \cite{K}, it follows that a toroidal embedding without self intersection is the same thing as a fine, saturated (f.s.) log smooth log scheme $(X,M_X)$ with log structure defined on the Zariski site of $X$. The ``characteristic monoid'' $\ov{M}_X:= M_X/\mathcal{O}^*_X$ is a constructible sheaf, constant on the strata of the associated stratification of $(X,U)$; its value at the generic point $\eta_Y$ of a stratum $Y$ is

\begin{align*}
\ov{M}_{X,\eta_Y} = (\sigma^{Y})^{\vee} \cap M^Y
\end{align*}

By basic results in the theory of log smoothness, or, as shown directly in \cite{KKMS}, a toroidal embedding admits \'etale locally an \'etale map to a toric variety, and, on a stratum $Y$, even a smooth map to $\Spec k[M^Y]$. 

By a morphism of toroidal embeddings $(X,U) \rightarrow (Y,V)$ we will simply mean a morphism of the associated log schemes; concretely, this means a morphism $X \rightarrow Y$ that takes the subsheaf $M_Y$ of $\c{O}_Y$ defined above into $M_X$.  The \emph{toroidal} morphisms of \cite{KKMS},\cite{AK} are the same as log smooth morphisms, defined by Kato in \cite{K}. We tend to prefer the latter terminology, as it has become more standard. The reader who prefers to avoid using logarithmic geometry can find a thorough treatment of toroidal morphisms in \cite[Section 1]{AK}. \\

We need the following observation, explained in \cite{KKMS}. We consider morphisms 
\begin{align*}
\lambda: \Spec k[[\NN]] \rightarrow X
\end{align*}
\noindent which take the generic point of $\Spec k[[\NN]]$ to $U$, and $\lambda(0)$ to $Y$. Here, $k[[\NN]]$ denotes the completion of $k[\NN] := k[t]$ with respect to its natural valuation, which we will denote by $\textup{ord}_0$. Then, for a divisor $D$ in $M^Y$, we get a pairing 
\begin{align*}
\langle \lambda,D \rangle = \textup{ord}_{0}\lambda^{*}D
\end{align*}
This way we obtain a map $\Hom{(\Spec k[[\NN]],X)} \rightarrow \sigma^Y$, which is actually surjective. So we may identify an integral point $v \in \sigma^Y \cap N^Y$ with an equivalence class of maps $\Spec k[[\NN]]$ to $X$, two maps being equivalent if and only if their order of intersection with each divisor is the same. We will abbreviate this equivalence class of maps by $v$ as well. Similarily, the cone $\sigma^Y$ itself can be identified with the image of $\Hom{(\Spec k[[\RR_+]], X)}$. Here, the ring $k[[\RR_+]]$ is the completion of $k[\RR_+]$ with respect to the valuation on $k[\RR_+]$ which takes a polynomial $\sum_{\alpha \in \RR_+} c^{\alpha}x^{\alpha}$ to $\textup{inf}\{\alpha: c_\alpha \ne 0\}$. The completion $k[[\RR_+]]$ is naturally a valuation ring as well; we denote its fraction field by $k[[\RR]]$. The relationship with $\sigma^Y$ is best explained through the following three observations. 

\begin{rem}
Suppose $V$ is an affine toric variety, corresponding to the cone $\sigma$ in the lattice $N$. Denote the dual lattice of $N$ by $M$ as usual, and denote by $\sigma^{\vee}$ the dual cone of $\sigma$, i.e $\{u \in M_\RR: \left \langle u,v \right \rangle \ge 0 \textup{ for all } v \in \sigma\}$, so that $V = \Spec k[ \sigma^{\vee} \cap M]$. An element $v \in \sigma \cap N$ is the same thing as a homomorphism of monoids $\sigma^{\vee} \cap M \rightarrow \NN$. We have 
\begin{align*}
\Hom_{\textup{Mon}} (\sigma^{\vee} \cap M,\NN) & = \Hom_{k-\textup{alg}} (k[\sigma^{\vee} \cap M], k[\NN]) = \\
& = \Hom_{\textup{Schemes}}(\Spec k[\NN] = \AA^1, V)
\end{align*} 
If $v$ in in the interior of $\sigma$, the image of $0$ under $\AA^1 \rightarrow V$ is precisely the torus fixed point of $V$. Composing $k[\sigma^{\vee} \cap M] \rightarrow k[\NN]$ with the completion $k[\NN] \rightarrow k[[\NN]]$ gives a morphism $\Spec k[[\NN]] \rightarrow V$ which defines precisely the same homomorphism $\sigma^{\vee} \cap M \rightarrow \NN$ as $v$. Thus, for an affine toric variety each equivalence class of morphisms $\Spec k[[\NN]] \rightarrow V$ has a canonical representative, obtained by completing the homomorphism $k[\NN] \rightarrow V$ correpsonding to $v$. 
\end{rem}
\begin{rem}
\label{rem:Rrepresentatives}
Suppose $V=V(\sigma)$ is the affine toric variety associated to the cone $\sigma$ in $N$, as in the preceeding remark. A similar description as the one given in the preceeding remark can in fact be given for any vector $v \in \sigma$ rather than just the integral ones. A vector $v$ in $\sigma$ (which we may well assume to be in the interior of $\sigma)$ corresponds to a homomorphism $\sigma^{\vee} \rightarrow \RR_+$, which induces (and is determined by) by restriction to a homomorphism $\sigma^{\vee} \cap M \rightarrow \RR_+$. This is the same data as a $k$-algebra homomorphism $k[\sigma^{\vee} \cap M] \rightarrow k[\RR_+]$ as above, which induces a map $k[\sigma^{\vee} \cap M] \rightarrow k[[\RR_+]]$, i.e a map $\lambda: \Spec k[[\RR_+]] \rightarrow V(\sigma)$. When $v$ is in the interior of $\sigma$, we have that under this morphism the closed point maps to the torus fixed point of $V(\sigma)$. Note that for any toric divisor $D$ of $V(\sigma)$, i.e any element $u \in M$, we have $\left \langle v,u \right \rangle = ord_{t=0} \lambda^{*} u$ by construction. We may thus identify vectors $v \in \sigma$ with morphisms 
\begin{align*}
\lambda: \Spec k[[\RR_+]] \rightarrow V(\sigma)
\end{align*}
such that $\lambda(\eta) \in \textup{ torus}$, and $\lambda(0) = \textup{torus fixed point}$, up to the equivalence relation that $ord_{t=0} \lambda^*$ induces the same homomorphism on $M$. We also have the analogue of the observation in the preceeding remark, that every equivalence class of a morphism $\Spec k[[\RR_+]] \rightarrow V(\sigma)$ has a unique representative obtained by completing $\Spec k[\RR_+] \rightarrow V(\sigma)$. 
\end{rem}
\begin{rem}
We can now combine the two remarks above with the fact that every toroidal embedding is \'etale locally (hence formally locally) isomorphic to a toric variety, to obtain that every interior vector $v \in \sigma^Y$, not necessarily integral, corresponds to an equivalence class of morphisms 
\begin{align*}
\lambda: \Spec k[[\RR_+]] \rightarrow V
\end{align*}
such that $\lambda(\eta) \in U$, $\lambda(0) \in Y$. Of course, there is no longer a canonical representative for a morphism in this equivalence class. 
\end{rem}

\begin{subsection}{Toroidal Stacks}
We would now like to transport the main points of the theory of toric stacks to toroidal embeddings. Though it is possible to give analogous definitions by working \'etale locally and modifying the appropriate results of \cite{GMtor}, we prefer not to work from scratch, and use the construction of the ``root stack" of Borne and Vistoli, \cite[Section 4]{BV}. We thus use the log scheme associated to a toroidal embedding.  

Given a log scheme $(X,M)$, the construction of Borne and Vistoli produces for any map of monoids $\ov{M} \rightarrow \ov{M}'$ which is \emph{Kummer}, i.e injective with finite cokernel, a log stack $(\mathcal{X},M')$ (which is tame, and Deligne-Mumford if the order of the cokernel is prime to the characteristic) mapping to $(X,M)$, with the morphism $M \rightarrow M'$ inducing the given map $\ov{M} \rightarrow \ov{M'}$. Furthremore, $(\mathcal{X},M')$ is the terminal object of the category of log schemes $(Y,N) \rightarrow (X,M)$ such that $M \rightarrow N$ factors through $M \rightarrow M'$. In the toroidal case, where we can choose \'etale locally on $X$ a chart $X \rightarrow \Spec k[\ov{M}]$, the stack $\mathcal{X}$ is described \'etale locally as a fiber product $X \times_{\Spec k[\ov{M}]} [\Spec k[\ov{M}']/K]$, where $K$ is the kernel of tori $\Spec k[\ov{M}^\gp] \rightarrow \Spec k[\ov{M}'^\gp]$. Put otherwise, if $\sigma$ is the dual monoid of $\ov{M}$ in $N=(\ov{M}^\gp)^\vee$, then $N \rightarrow N' = (\ov{M'}^\gp)^\vee$ is a finite index inclusion, and $\mathcal{X}$ is \'etale locally isomorphic to $X \times_{\AA(\sigma,N)} \AA(\sigma,N,N')$.    

Starting from the log scheme $(X,M_X)$ associated to a toroidal embedding, and a Kummer extension $\ov{M}_X \rightarrow \ov{M'}_X$, the stack $(\mathcal{X},M'_\mathcal{X})$ is log smooth with Zariski log structure. We will refer to it as a toroidal stack (without self intersection). To translate the Borne-Vistoli formalism to a formalism more analogous to the toric formalism of the previous section, given a toroidal embedding $(X,U)$, we will assign for each cone $\sigma^{Y}$ in $C(X)$ a sublattice $N_{\sigma^Y} \subset N^Y \cap \Span(\sigma^Y)$ which is injective with finite cokernel, and with the property that $N_{\tau} = \Span \tau \cap N_{\sigma}$ for a face $\tau$ of $\sigma$. Recall that the sheaf $\ov{M}_X$ is constructible, constant on the strata of $(X,U)$, with value $M^Y$ on $Y$; therefore, dualizing the inclusions $N_{\sigma^Y} \rightarrow N^Y$ gives exactly a Kummer extension $\ov{M}'_X$ of $\ov{M}_X$. Thus, a compatible triple $(C(X),\{N^Y\},\{N_{\sigma^Y}\})$, as $Y$ ranges through the strata of $(X,U)$  produces a toroidal stack over $X$ by the \cite{BV} construction.

Suppose now $C(X')$ is a subdivision of $C(X)$. The process outlined in \cite{KKMS} produces a log blowup $X' \rightarrow X$ with cone complex $C(X')$. Combining such subdivisions with Kummer extensions produces triples $(C(X'),\{N^Z\},\{N_{\sigma^Z}\})$, with $Z$ running through the strata of $(X',U)$, and this data yields a stack $\c{X'}$ mapping to $X$. We call such a triple a stacky subdivision of $C(X)$, and refer to it as the cone complex of $\c{X'}$. The notation $(C(X'),\{N^Z\},\{N_{\sigma^Z}\})$ may be cumbersome, but we find it very useful, as doing geometry on this stacky polyhedral complex is much more intuitive than doing geometry on the dual sheaf of monoids. The following lemma is a combination of \cite{KKMS} and \cite{BV}. 

\begin{lemma}
Let $X$ be a toroidal embedding, $\c{X'}$ the log stack corresponding to a triple $(C(X'),\{N^Z\},\{N_{\sigma^Z}\})$ which is a stacky subdivision of $C(X)$. Let $T \rightarrow X$ be a map from a toroidal embedding $T$, and assume that $C(T) \rightarrow C(X)$ factors through $(C(X'),\{N^Z\},\{N_{\sigma^Z}\})$. Then there exists a unique lift of $T \rightarrow X$ to $T \rightarrow \c{X'}$ which induces the given map $C(T) \rightarrow (C(X'),\{N^Z\},\{N_{\sigma^Z}\})$.
\end{lemma} 

\begin{proof}
Note first that by \cite{KKMS}, the map $T \rightarrow X$ uniquely factors through $X'$ where $X'$ is the toroidal embedding with cone complex $C(X')$ constructed in \cite{KKMS}. Then $\c{X'} \rightarrow X'$ is the log stack associated to the Kummer extension $(C(X'),N^Z,N_{\sigma^Z})$ of $C(X')$ and the result follows from \cite{BV}. 
\end{proof}

The combinatorial definition of weak semistability \ref{def: weaklysemistable} for toric fans does not make use of the global lattices in which the fans live, so it applies without change to maps of cone complexes. The connection with algebraic geometry is as follows: 

\begin{lemma}
Let $f: X \rightarrow S$ be a log smooth morphism of toroidal embeddings, and assume that the map of cone complexes $C(X) \rightarrow C(S)$ is weakly semistable. Then $X \rightarrow S$ is flat with reduced fibers. 
\end{lemma}

\begin{proof}
For the proof, we use some basic results from the theory of logarithmic geometry. It suffices to check that the conditions hold around every point $x \in X$; as the question is \'etale local on $X$, we may assume that $X$ and $S$ have a single closed stratum, and that $x$ and $f(x)$ belong to the closed strata. In terms of cone complexes this means that $C(X)$ and $C(S)$ consist of a single cone, and the assumption of weak semistability translates to the fact that the map $\sigma^X \cap N^X \rightarrow \sigma^S \cap N^S$ is surjective; this in turns implies that the dual map $\overline{M}_{S,f(x)}^\gp \rightarrow \overline{M}_{X,x}^\gp$ is injective with torsion-free cokernel. Then, by replacing $X$ with an \'etale cover, we can choose a neat chart for $f$, which we can even take to be a characteristic chart (see for instance by \cite[Theorem 2.4.4]{oguslog}); then, we have a diagram
\begin{align*}
\begin{xymatrix}
{X \ar[r] \ar[d] & \Spec \ZZ[\overline{M}_{X,x}] \ar[d] \\  S \ar[r] & \Spec \ZZ[\overline{M}_{S,f(x)}]} 
\end{xymatrix}
\end{align*}
From the toric case, we know that $\Spec \ZZ[\overline{M}_{X,x}] \rightarrow \Spec \ZZ[\overline{M}_{S,f(x)}]$ is flat and has reduced fibers, and, as the chart is neat, the map $X \rightarrow S \times_{\Spec \ZZ[\overline{M}_{S,f(x)}]} \Spec \ZZ[\overline{M}_{X,x}]$ is smooth. Therefore, the composition $X \rightarrow S$ is also flat with reduced fibers. 
\end{proof}

 By the fiber product of morphisms of toroidal embeddings we will mean the fiber product \emph{in the category of} f.s. log schemes. We will denote this by $(X \times_S T)_{\textup{tor}}$ to avoid confusion, and write $X \times_S T$ for the schematic fiber product. For log smooth morphisms, this fiber product is a toroidal embedding, and the underlying scheme of this fiber product is connected to the fiber product of the underlying schemes in the way explained in \ref{subsection: fiberproducts} -- it is the normalization of the closure of the main component of the underlying fiber product. This follows from the existence of the \'etale local maps to toric varieties, as the statement can be checked \'etale locally. A consequence of weak semistability is therefore that the fiber product in the category of f.s. log schemes has the same underlying scheme as the fiber product of the underlying schemes. The cone complex of the fiber product $(X \times_S T)_{\textup{tor}}$ is the fiber product of cone complexes $C(X) \times_{C(S)} C(T)$. This again follows from the local case.  
\end{subsection}

\begin{subsection}{The Main Construction}

The toric construction expalined in section $2$ carries over to the toroidal case with minimal changes, by replacing the fans of $X$ and $S$ with the cone complexes $C(X),C(S)$. The morphism $X \rightarrow S$ induces a morphism $p:C(X) \rightarrow C(S)$ by composing a map $k[[ \RR_+]] \rightarrow X$ with $X \rightarrow S$. So suppose a cone $\kappa \in C(S)$, and a point $w \in \kappa$ are given. We consider 
\begin{align*}
N_0(w) = \{\sigma: \exists ! v \in \sigma^o \textup{ such that } p(v) = w\}
\end{align*}

\begin{lemma}
\label{lem:convex}
The cones $\{w: N_0(w) = \textup{ fixed} \}$ are strictly convex, and form a subdivision of $\kappa$. 
\end{lemma}

\begin{proof}
We adapt the proof of the toric case \ref{theorem:maintoric}. The question is local on $S$, so we may assume that $C(S)$ is a single cone $\kappa$. We take two vectors $w,w'$ in $\kappa$ for which $N_0(w)=N_0(w')$, and try to show that we have $N_0(tw+(1-t)w')=N_0(w)$ for all $t \in [0,1]$. As above, we may assume that this condition fails for some $t$ and derive a contradiction, and even take for $t$ the minimal element of $[0,1]$ for which the condition fails. So we may replace $\kappa$ by the interval $[w,w']$, and $C(X)$ by its fiber over $[w,w']$, which we will denote by $C(X)$ as well for simplicity. Put $w''=tw+(1-t)w'$. As above, we can choose an element in $N_0(w'')-N_0(w)$, which corresponds to a vertex $v''$ in $C(X)$. The key step in the toric proof is that the star of $v''$ in $C(X)$ intersects the fiber of $C(X)$ over $[w,w'')$, which follows from the properness of the map. This statement is not immediately clear in the toroidal situation, but we claim it is nevertheless still correct. To see this, pick a family of maps 
\begin{align*}
\Spec k[[\RR_+]] \times [0,1] \rightarrow S
\end{align*}
corresponding to the interval $[w,w']$ in $\kappa$, and a lift 
\begin{align*}
\Spec k[[\RR_+]] \rightarrow X
\end{align*}
corresponding to $v''$. We abusively denote the lift by $v''$ as well. Let $x=v''(0) \in X$. Since $X \rightarrow S$ is log smooth, we may choose a chart 
\begin{align*}
\xymatrix{ (X,x) \ar[r]^f \ar[d]_p & (V,v) \ar[d]^\pi \\ (S,p(x)) \ar[r]_g & (W,\pi(v))}
\end{align*} 
\noindent for the morphism $p$, where: the horizontal morphisms $f,g$ are \'etale; $V,W$ are toric varieties, $\pi$ is a toric morphism, and $v,\pi(v)$ are special points in the torus orbits; and the morphism $N_\RR \rightarrow Q_\RR$ is surjective, where $N,Q$ are the lattices of $V$ and $W$ respectively. Let $[z,z']$ denote the interval corresponding to $[w,w']$ in $Q_\RR$ under $g$, and let $y'' \in N_\RR$ denote the element corresponding to $v''$ under $f$. Since $N_\RR \rightarrow Q_\RR$ is surjective, we may lift $[z,z']$ to an interval $[y,y'] \in N_\RR$ lying over $[z,z']$, with $y'' \in [y,y']$. In other words, if we denote by $T(V),T(W)$ the tori of $V$ and $W$ respectively, we get that the family
\begin{align*}
\Spec k[[\RR]] \times [0,1] \rightarrow  T(W) \subset W
\end{align*}
\noindent corresponding to $[z,z']$ lifts to a family of maps 
\begin{align*}
\Spec k[[\RR]] \times [0,1] \rightarrow  T(V) \subset V
\end{align*}
\noindent which under $\pi$ projects to $[z,z']$, and with $y''$ the map over $z''$, i.e such that the morphism at $t \in [0,1]$ is $y''$.  Composing with the inverse of the isomorphism $\mathcal{\hat{O}}_{X,x} \cong \mathcal{\hat{O}}_{V,v}$ we get a family of maps 
\begin{align*}
\Spec k[[\RR]] \times [0,1] \rightarrow X
\end{align*}
\noindent which under $p$ compose to the original family $\Spec k[[\RR]] \rightarrow S$ corresponding to $[w,w']$. Now, for each $s \in [0,1]$, the map $\Spec k[[\RR]] \rightarrow X$ extends to a map $\Spec k[[\RR_+]] \rightarrow X$ by the properness of $X \rightarrow S$. But such a map corresponds to an element of $C(X)$, so we get a family of elements in $C(X)$ which at $s=t$ specialize to $v''$. Thus, the star of $v''$ in $C(X)$ contains a lift of the line segment $[w,w'')$, and the same argument as in the toric case goes through. 
\end{proof}

Theorem \ref{theorem:maintoric} now carries through without any change in the proof, once we consider the appropriate generalization of the category $\mathcal{C}$ in the toroidal setting. We fix a proper, surjective log smooth morphism $X \rightarrow S$, which gives a cone complex morphism $C(X) \rightarrow C(S)$. We consider the subdivision of $S$ determined by the cones $\{w \in C(S): N_0(w) = \textup{ constant}\}$, and the subdivision of $C(X)$ given by $p^{-1}(\{w \in C(S) :N_0(w) = \textup{ constant}\}) \cap \sigma$. For a cone $\kappa$ whose interior is given by the collection $\{w: N_0(w)=\{\sigma_i\}_1^n\}$, we take the sublattice $Q_\kappa'$of $Q_\kappa$ to be the lattice generated by the elements in $\cap_{i=1}^{n} p(\sigma_i \cap N_{\sigma_i})$, and for a cone $\sigma_i \cap p^{-1}(\kappa)$ we take the sublattice $N'_{\sigma_i \cap p^{-1}(\kappa)} = p^{-1}(Q'_\kappa)$. This construction produces a log smooth morphism of toroidal stacks $\mathcal{X} \rightarrow \mathcal{S}$. \\

\begin{defn}
Let $\mathcal{C}_t$ be the category whose objects are diagrams

\begin{align*}
 \xymatrix{Y \ar[r]^j \ar[d]_\pi & X \ar[d]^p \\ T \ar[r]_i & S} 
\end{align*}

such that 
\begin{itemize}
\item $Y,T$ are toroidal embeddings and $j,\pi,i$ are log smooth morphisms. \\
\item The map $i$ is an alteration. \\
\item $Y$ is a modification of the fiber product $(X \times_S T)_{\textup{tor}}$. \\ 
\item $\pi$ is weakly semistable. \\
\end{itemize}
A morphism in $\mathcal{C}_t$ is a commutative diagram 
\begin{align*}
\xymatrix{Y' \ar[r] \ar[d] & Y \ar[d] \\ T '\ar[r] & T}
\end{align*}
\noindent which commutes with the morphisms to $p:X \rightarrow S$.
\end{defn} 

Then we have 

\begin{theorem}
The family $\mathcal{X} \rightarrow \mathcal{S}$ is the terminal object of $C_t$. 
\end{theorem}

\begin{proof}
The proof of the toric case in \ref{theorem:maintoric} carries through in exactly the same way. 
\end{proof}

As a corollary, we get

\newtheorem*{theorem:main}{Theorem \ref{theorem:main}}
\begin{theorem:main}[Universal Weak Semistable Reduction]
Let $X \rightarrow S$ be any proper, surjective, log smooth morphism of toroidal embeddings. Then, there exists a commutative diagram 
\begin{align*}
\xymatrix{ \mathcal{X} \ar[r] \ar[d] & X \ar[d] \\ \mathcal{S} \ar[r] & S } 
\end{align*}
where $\mathcal{X} \rightarrow \mathcal{S}$ is a representable morphism of tame algebraic stacks, such that for any diagram
\begin{align*}
\xymatrix{ Y \ar[r] \ar[d] & X \ar[d] \\ T \ar[r] & S}
\end{align*}
where $T \rightarrow S$ a toroidal alteration; $Y$ is a modification of the fiber product $(X \times_S T)_{\textup{tor}}$; and $Y \rightarrow T$ is weakly semistable, the morphism $Y \rightarrow T$ factors uniquely through $\mathcal{X} \rightarrow \mathcal{S}$. Furthermore, $\mathcal{X} \times_\mathcal{S} T \rightarrow T$ is weakly semistable.  
\end{theorem:main}      

\end{subsection}

\begin{subsection}{Generalizations}
We note first of all that the assumption that $T \rightarrow S$ is a toroidal alteration is not important. The theorem holds for an arbitrary log map $T \rightarrow S$, with the same proof. The reason to state the theorem with $T \rightarrow S$ an alteration is to emphasize that $\mathcal{S} \rightarrow S$ is terminal in the category of alterations so as to make contact with the statement of \cite{AK}. The assumption that $k$ is algebraically closed is also not needed, with the caveat that in the non-algebraically closed case we work with toric varieties with a split torus. \\

The assumption that $X$ and $S$ are toroidal is more serious, though it is important only in the following ways: (1) A toroidal embedding has a ``fan'', i.e., its cone complex $C(X)$; (2) a map $X \rightarrow Y$ induces a map of fans $C(X) \rightarrow C(Y)$; (3) weak semistability of $C(X) \rightarrow C(S)$ implies weak semistability of $X \rightarrow S$. In short, for toroidal embeddings there is a functorial theory of fans, and the fan captures enough local information for $X \rightarrow S$. In general, there are many constructions of fans of logarithmic schemes of varying levels of generality -- for instance, in \cite{KKMS},\cite{Ktor}, or \cite{ACMUW}; all constructions are essentially equivalent, and they all satisfy (3) for log smooth morphisms; the constructions however are not functorial. It is unclear for what classes of morphisms a theory of fans satisfying (1),(2),(3) can exist. \\

Since we could not see any reasonable generalization to a wider class of log schemes, we decided to keep the statements in the language of toroidal embeddings as much as possible, as it is more widely known. However, it may be worth noting that given any context for which assumptions (1), (2), and (3) do hold, the main theorem \ref{theorem:main} also holds. For instance, the theorem holds for log curves over a log point. Starting with a log curve $X \rightarrow S$ over a log point, one often needs to perform a log blowup $Y \rightarrow X$; this produces a curve $Y \rightarrow S$, but this is not a log curve (in the sense of F.Kato \cite{fK}), as $Y \rightarrow S$ is not weakly semistable. However, the construction still applies and gives a universal log curve $\c{Y} \rightarrow \c{S}$ over a stacky blowup of the base.   

\end{subsection}
\end{section}
\bibliography{refs}{}
\end{document}